\documentstyle{amsppt}

\magnification=\magstep1
\hcorrection{.25in}
\advance\voffset.25in
\advance\vsize-.5in

\define\e{\varepsilon}

\define\th{\theta}

\define\conv{\text{conv}}
\define\etc{, \dots ,}
\define \sumi{\sum_{i=1}^M}

\define\rn{$\Bbb R^n \,$}
\define\Rm{$\Bbb R^M \,$}
\define\RM{\Bbb R^M}
\define\Rn{\Bbb R^n}

\define\BM{B^M_2}
\define\nor #1{\left \| #1 \right \|}
\define\enor #1{\Bbb E \, \nor{#1}}

\topmatter

\title
Almost orthogonal submatrices of an orthogonal matrix
\endtitle

\author
M. Rudelson
\endauthor

\thanks
Research supported in part by a grant of the US--Israel BSF.
Research at MSRI is supported in part  by NSF grant DMS-9022140.
\endthanks


\address
\hskip-\parindent
Institute of Mathematics, 
The Hebrew University of Jerusalem, \hfill \break
Givat Ram, 91904  Jerusalem, Israel \hfill \break
and Mathematical Sciences Research Institute, \hfill \break
1000 Centennial Drive, Berkeley, CA 94720, USA
\endaddress

\email
mark\@math.huji.ac.il, \ mark\@msri.org
\endemail

\abstract
Let  $t \ge 1$ and let $n, \ M$ be natural numbers, $n < M$.
Let $A=(a_{i,j})$ be an $n \times M$ matrix  whose rows are orthonormal.
Suppose that for all $j$
$$
\sqrt{\frac{M}{n}} \cdot \left ( \sum_{i=1}^n a_{i,j}^2 \right )^{1/2} \le t . 
$$
Using majorizing measure estimates we prove that
for every $\e >0$ there exists a set $I \subset \{ 1, \dots, M \}$
of cardinality at most
$$
 C \cdot \frac{t^2}{\e^2} \cdot n \cdot \log n
$$
 so that for all $x \in \ell_2^n$  
$$
(1-\e ) \cdot \nor{x} \le \sqrt{\frac{M}{|I|}} \cdot \nor{R_I A^T x}
\le (1+ \e ) \cdot \nor{x}. 
$$
Here $R_I: \Bbb R^M \to \Bbb R^M$ is the orthogonal projection onto the
space span$\{ e_i \ | \ i \in I \}$, where $\{ e_i \}_{i=1}^M$ is the
standard basis of $\ell_2^M$.

\endabstract

\endtopmatter

\document

\head 1. Introduction \endhead
We consider the following problem, posed by B. Kashin and L. Tzafriri 
\cite{K-T}:
\flushpar {\sl
Let  $\e > 0$ and let $n, \ M$ be natural numbers, $n < M$.
Given an $n \times M$ matrix $A$ whose rows are orthonormal,
what is the smallest cardinality  $L(A,\e)$ of a subset 
$I \subset \{ 1, \dots, M \}$  so that for all $x \in \ell_2^n$  }
$$
(1-\e ) \cdot \nor{x} \le \sqrt{\frac{M}{|I|}} \cdot \nor{R_I A^T x}
\le (1+ \e ) \cdot \nor{x}. \tag 1.1
$$
Here $R_I: \Bbb R^M \to \Bbb R^M$ is the orthogonal projection onto the
space span$\{ e_i \ | \ i \in I \}$, where $\{ e_i \}_{i=1}^M$ is the
standard basis of $\Bbb R^M$. 
Throughout this paper we denote by $\nor{\cdot}$ the standard $\ell_2$-norm
and by $|I|$ the cardinality of a set $I$.

Under an additional assumption that all the entries of $A$ have the same
 absolute value $1/\sqrt{M}$    Kashin and  Tzafriri
proved that 
$$
L(A,\e) \le \frac{c}{\e^4} \cdot n^2 \log n .  \tag 1.2
$$
Moreover, their proof shows that a random subset $I$ of this cardinality
satisfies (1.1) with probability close to 1. 
Clearly, the estimate (1.2) is not optimal.
The example of random selection of columns of a rectangular Walsh matrix,
considered by Kashin and  Tzafriri suggests that the possible upper bound 
could be
$$
L(A,\e) \le C(\e) \cdot n \log n .  \tag 1.3
$$
From the other side, simple examples (\cite{K-T}, \cite{R}) show that 
the estimate (1.3) is the best one can obtain by the random selection method.

As it was mentioned in \cite{R}, the Kashin and Tzafriri problem is dual to 
that of finding an approximate John's decomposition. 
Entropy estimates used in \cite{R} for the last problem 
enabled to improve (1.2).
More precisely, let $t \ge 1$ and suppose that the matrix $A$ satisfies
$$
\sqrt{\frac{M}{n}} \cdot \left ( \sum_{i=1}^n a_{i,j}^2 \right )^{1/2} \le t. 
$$
for all $j= 1 \etc M$.
Then 
$$
L(A,\e) \le C(\e ) \cdot t^2 \cdot n \log^3 n.
$$ 

In order to improve this estimate one can use majorizing measures instead
of entropy estimates.
The method of majorizing measures, developed by Talagrand (\cite{L-T},
\cite{T1}), is extremely useful in obtaining estimates of stochastic 
processes, related to random selection.
A random process, similar to that arising in the Kashin and Tzafriri problem 
was considered by Talagrand \cite{T2} for the problem of embedding of
a finite dimensional subspace of $L_p$ into $\ell_p^N$. 
For this kind of processes Talagrand introduced a special method of 
constructing majorizing measures.
This method ($s$-separated trees) can be used to prove an estimate
$$
L(A,\e) \le C(\e ) \cdot t^2 \cdot n \log n \cdot (\log \log n)^2
$$ 
for the Kashin and Tzafriri problem.
It is unlikely that the $(\log \log n)^2$ factor can  be removed by 
a modification of the $s$-separated trees method.
However, using a different approach based on the explicit construction of a 
partition tree, we obtained a sharper
estimate.
More precisely, we prove the following
\proclaim{Theorem}
Let $t \ge 1$ and let $A=(a_{i,j})$ be an $n \times M$
matrix, whose rows are orthonormal.
Suppose that for all $j$
$$
\sqrt{\frac{M}{n}} \cdot \left ( \sum_{i=1}^n a_{i,j}^2 \right )^{1/2} \le t. 
 \tag 1.4
$$
Then for every $\e >0$ there exists a set $I \subset \{ 1, \dots, M \}$
so that 
$$
|I| \le C \cdot \frac{t^2}{\e^2} \cdot n \cdot \log n  \tag 1.5
$$
and for all $x \in \Rn$
$$
(1-\e ) \cdot \nor{x} \le \sqrt{\frac{M}{|I|}} \cdot \nor{R_I A^T x}
\le (1+ \e ) \cdot \nor{x}.  \tag 1.6
$$
\endproclaim
Throughout this paper $C,c$ etc. denote absolute constants whose value may
 change from line to line.

The main part of the proof is the proof of Lemma 1 below.
Our original proof of this lemma used the direct construction
of the majorizing measure.
It included an explicit construction of a sequence of partitions and 
putting  weights on the elements of each partition.
This scheme is based on the Talagrand and Zinn's proof of the 
majorizing measure theorem of Fernique (Proposition 2.3 and Theorem 2.5 
\cite{T4}).
The proof was rather involved, since we had to approximate the natural
metric of a random process by a family of metrics depending on the
elements of the partition.
After we had shown our proof to  M.~Talagrand, he pointed out that the 
explicit construction of the partition tree may be substituted by
applying his general majorizing measure construction (Theorems 4.2, 4.3 
and Proposition 4.4 \cite{T4}).
This resulted in a considerable simplification of the proof.
We present here the argument suggested by Talagrand.

By the duality between the Kashin and Tzafriri problem and approximate John's 
decompositions, we have the following
\proclaim{Corollary}
Let B be a convex body in \rn and let $\e>0$.
There exists a convex body $K \subset \Rn$, so that ${\text d}(K,B) \le 1+\e$
and the number of contact points of $K$ with its John ellipsoid is less
than
$$ 
m(n,\e) = C(\e) \cdot n \cdot \log  n. 
$$
\endproclaim

\head 2. The random selection method \endhead
Clearly, we may assume that 
$M \ge C \cdot \frac{t^2}{\e^2} \cdot n \cdot \log n $ for some 
absolute constant $C$.

The proof of the Theorem is based on the following iteration procedure.
Let $A=(a_{i,j})$ be an $ n \times M$ matrix , satisfying (1.4),
We define a sequence $\{ \e_i \}_{i=1}^M$ of independent Bernoulli
variables taking values $\pm 1$ with probability $1/2$ and put
$$
I_1 = \{ i \bigm | \e_i=1 \}.
$$
Then
$$
\frac{M}{2} \cdot \left ( 1- \frac{1}{\sqrt{M}} \right ) \le |I_1|  \le 
\frac{M}{2}  \tag 2.1
$$
with probability at least $1/4$.
Define 
$$
W=A \Bbb R^n
$$
and denote by $w(1) \etc w(M)$ the coordinates of a vector $w$.
We have to estimate
$$
\align
\sup_{x \in B_2^n} \left | 2 \nor{R_I A x}^2 - \nor{x}^2 \right | &=
\sup_{w \in W \cap \BM} \left | 2 \cdot \sum_{i \in I} w^2(i) - \sumi w^2(i)
\right | \\
&= \sup_{w \in W \cap \BM} \left | \sumi \e_i  w^2(i) \right |.
\endalign
$$
Denote by $\Bbb E X$ the expectation of a random variable $X$.
The key step of the proof is the following
\proclaim{Lemma 1}
Let $W$ be an $n$-dimensional subspace of \Rm. 
Let $\e_1 \etc \e_M$ be independent Bernoulli variables taking values 
$\pm 1$ with probability $1/2$.
Then 
$$
\Bbb E \ \sup_{w \in W \cap \BM} \left | \sumi \e_i w^2(i) \right | 
\le C \sqrt{\log M} \cdot
\nor{P_W: \ell_1^M \to \ell_2^M}.
$$
Here $P_W: \RM \to \RM$ is the orthogonal projection onto $W$.
\endproclaim 

From (1.4) it follows that
$$
\nor{P_W: \ell_1^M \to \ell_2^M} \le t \cdot \sqrt{\frac{n}{M}},
$$
so by Lemma 1 and Chebychev's inequality we have
$$
\sup_{x \in B_2^n} \left | 2 \nor{R_{I_1} A x}^2 - \nor{x}^2 \right | \le
C \cdot t \cdot \sqrt{\frac{n}{M}} \cdot \sqrt{\log M}  \tag 2.2
$$
with probability more than $3/4$.
Thus, there exists a set $I_1 \subset \{ 1 \etc M \}$ satisfying (2.1) and
(2.2).

Repeating this procedure, we obtain a sequence of sets
$\{ 1 \etc M \} = I_0 \supset I_1 \supset \ldots \supset I_s$ so that
$$
\frac{|I_k|}{2} \cdot \left ( 1- \frac{1}{\sqrt{|I_k|}} \right ) \le
|I_{k+1}| \le \frac{|I_k|}{2} \tag 2.3
$$
and
$$
\sup_{x \in B_2^n} \left ( 2^k \nor{R_{I_k} A x}^2 - 
 2^{k-1} \nor{R_{I_{k-1}} A x}^2 \right ) \le
C \cdot t \cdot \sqrt{\frac{n}{M/2^k}} \cdot \sqrt{\log |I_{k-1}|}. \tag 2.4
$$

Indeed, at each step of induction we have
$$
\frac12 \nor{x} \le 2^{\frac{k-1}{2}} \nor{R_{I_{k-1}} A x} \le 
\frac32 \nor{x}. \tag 2.5
$$
Assume for simplicity that $I_{k-1}= \{ 1 \etc m \}$ for some $m<M$.
Let $W_k = R_{I_{k-1}} A \RM \subset \Bbb R^m$ and let $
P_{W_k}: \Bbb R^m \to \Bbb R^m$ be the orthogonal projection onto $W_k$.
Then 
$$
2^{\frac{k-1}{2}} R_{I_{k-1}} A B_2^n \subset \frac32  B_2^m \cap W_k,
$$
so for a random set $I_k \subset \{ 1 \etc m \}$ we have
$$
\align
\Bbb E \, \sup_{x \in B_2^n} \left ( 2^k \nor{R_{I_k} A x}^2 - 
 2^{k-1} \nor{R_{I_{k-1}} A x}^2 \right ) \le 
&\Bbb E \, \sup_{w \in \frac32  B_2^m \cap W_k} \sum_{i=1}^m \e_i w^2(i) 
\le \\
\frac94 \Bbb E \, \sup_{w \in B_2^m \cap W_k} \sum_{i=1}^m \e_i w^2(i).
\endalign
$$
To apply Lemma 1 we need to compute $\nor{P_{W_k}: \ell_1^m \to \ell_2^m}$.
By (2.5) we have 
$$ 
\nor{P_{W_k}: \ell_1^m \to \ell_2^m} \le 
2 \cdot \nor{\left ( 2^{\frac{k-1}{2}} R_{I_{k-1}} A \right )^*: 
\ell_1^M \to \ell_2^n} 
\le 2^{\frac{k+1}{2}} \cdot t \cdot \sqrt{\frac{n}{M}}.
$$
Now (2.4) follows from Lemma 1 and Chebychev's inequality.

Summing up inequalities (2.4) we get 
$$
\aligned
&\sup_{x \in B_2^n} \left | 2^s \nor{R_{I_s} A x}^2 - \nor{x}^2 \right | \le
C \cdot t \cdot \sqrt{\frac{n}{M/2^s}} \cdot \sqrt{\log |I_s|} \le \\
& C \cdot t \cdot \sqrt{\frac{n}{M/2^s}} \cdot \sqrt{\log \frac{M}{2^s}}. 
\endaligned \tag 2.6
$$
We proceed until the last expression is greater than $\e /2$.
In this case
$$
c \cdot \frac{t^2}{\e^2} \cdot n \cdot \log n \le \frac{M}{2^s} \le
C \cdot \frac{t^2}{\e^2} \cdot n \cdot \log n .
$$
From (2.3) it follows  that
$$
\frac{M}{2^s} \cdot \left ( 1- \frac{4}{\sqrt{|I_s|}} \right )  \le
|I_s| \le \frac{M}{2^s}, 
$$
so we obtain (1.5) and
$$
\frac{M}{|I_s|} \cdot \left ( 1- 
\left (c \cdot \frac{t^2}{\e^2} \cdot n \cdot \log n \right )^{-1/2}
 \right )  \le
2^s \le \frac{M}{|I_s|}.
$$
Then, (2.6) implies that 
$$
\sup_{x \in B_2^n} 
\left | \frac{M}{|I_s|} \cdot  \nor{R_{I_s} A x}^2 - \nor{x}^2 \right | \le
\e
$$
and this completes the proof of the Theorem. \hfil \qed
\demo{Remark}
The random selection method was used first by Talagrand \cite{T3} to simplify
the construction of embedding of a finite dimensional subspace of $L_1$
into $\ell_1^N$.
The original construction of Bourgain, Lindenstrauss and Milman used
the empirical distribution method instead of it.
The advantage of the random selection is that it enables to deal with
random processes having a subgaussian tail estimate, rather than with 
general Bernoulli processes.
\enddemo

\head 3. Construction of the majorizing measure. \endhead
The proof of  Lemma 1 uses the majorizing measure theorem of Talagrand
\cite{T1}, \cite{T4}.
This theorem provides a bound to
$$
\Bbb E \sup_{t \in T} X_t
$$
for a subgaussian process $X_t$ indexed by points of a metic space $T$ with
a metric $d$ through the geometry of this space.
However it turns out that the space $T$ does not have to be assumed metric.
The same proof works in the case when $d$ is a quasimetric, i.e. if there
exists a constant A such that for any $t, \bar t, s \in T$ 
$$
d(t, \bar t) \le A \cdot \bigl (d(t,s) + d(s, \bar t) \bigr ).
$$
We use the following version of 
\proclaim{Majorizing measure theorem}
Let $(T,d)$ be a quasimetric space.
Let $(X_t )_{t \in T}$ be a collection of mean 0 random variables with the
subgaussian tail estimate
$$
 \Cal P \, \{| X_t - X_{\bar t}| > a \} \le 
\exp \left ( -c \frac{a^2}{d^2 (t, \bar t)} \right ),
$$
for all $a>0$.
Let $r>1$ and let $k_0$ be a natural number so that the diameter of $T$ is 
less than $r^{-k_0}$.
Let $\{\varphi_k \}_{k=k_0}^{\infty}$ be a sequence of functions from 
$T$ to $\Bbb R^+$,  uniformly bounded by a constant depending only on $r$.
Assume that there exists $\sigma >0$ so that for any $k$ the functions 
$\varphi_k$ satisfy the following condition: \hfil \break
\par\flushpar
 for any $s \in T$ and for any points 
$t_1 \etc t_N \in B_{r^{-k}}(s)$ with mutual distances at least $r^{-k-1}$ 
one has
$$
\max_{j=1 \etc N} \varphi_{k+2}(t_j) \ge \varphi_k(s) + \sigma \cdot r^{-k}
\cdot \sqrt{\log N}. \tag 3.1
$$
Then 
$$
\Bbb E \sup_{t \in T} X_t \le C(r) \cdot \sigma^{-1}.
$$
\endproclaim

This version may be obtained as a combination of the majorizing
measure theorem of Fernique \cite{L-T} and the general majorizing
measure construction of Talagrand (Theorems 2.1 and 2.2 \cite{T1} or
Theorems 4.2, 4.3 and Proposition 4.4 \cite{T4}).
 
To prove Lemma 1 we need some estimates of covering numbers.
Denote by $N(B,d,\e)$ the $\e$-entropy of $B$, i.e. the number of 
$\e$-balls in the  (quasi--) metric $d$ needed to cover the body $B$.
We use the following 
\proclaim{Lemma 2}
Let $W$ be an $n$-dimensional subspace of \Rm and let $P_W$ be the 
orthogonal projection onto $W$.
\roster
\item $\e \sqrt{\log N(\BM \cap W, \nor{\cdot}_{\infty}, \e)} \le
C \cdot \nor{P_W: \ell_1^M \to \ell_2^M} \cdot \sqrt{\log M}$;
\item Let $\nor{\cdot}_{\Cal E}$ be a norm defined by
$$ 
\nor{x}_{\Cal E}= \left ( \sumi x^2(i) \cdot a_i^2 \right )^{1/2}.
$$
Then 
$$
\e \sqrt{\log N(\BM \cap W, \nor{\cdot}_{\Cal E}, \e)} \le
C \cdot 
\nor{P_W: \ell_1^M \to \ell_2^M} \cdot \left ( \sumi a_i^2 \right )^{1/2}.
$$
\endroster 
\endproclaim
\demo{Proof}
Both statements follow from the dual Sudakov minoration \cite{L-T}.

(1) Let $g$ be the standard Gaussian vector in \Rm.
Then $P_W g$ is the standard Gaussian vector in the space $W$.
So,
$$
\align
&\e \sqrt{\log N(\BM \cap W, \nor{\cdot}_{\infty}, \e)} \le C \cdot
\Bbb E \, \nor{P_W g}_{\infty} = C \cdot \Bbb E \, \max_{j=1 \etc M} 
| \langle P_W g, e_j \rangle | \le \\
& C \cdot \sqrt{\log M} \cdot \max_{j=1 \etc M} \nor{P_W e_j} =
C \cdot \sqrt{\log M} \cdot \nor{P_W: \ell_1^M \to \ell_2^M}. \qed
\endalign 
$$

(2) Again dual Sudakov minoration gives
$$
\align
&\e \sqrt{\log N(\BM \cap W, \nor{\cdot}_{\Cal E}, \e)} \le C \cdot
\Bbb E \, \nor{P_W g}_{\Cal E} \le C \cdot 
\left ( \Bbb E \, \nor{P_W g}_{\Cal E}^2 \right )^{1/2} = \\
& C \cdot 
\left ( \Bbb E \, \sumi \langle P_W g, e_i \rangle^2 a_i^2\right )^{1/2}
\le C \cdot \max_{i=1 \etc M} \nor{P_W e_i} \cdot 
 \left ( \sumi a_i^2 \right )^{1/2}. \qed
\endalign
$$
\enddemo

\demo{Proof of Lemma 1}
Denote 
$$
W_1=\BM \cap W.
$$
We have to estimate the expectation of the supremum over all $w \in W_1$
of a random process 
$$
V_w = \sumi \e_i w^2(i).
$$
The process $V_w$ has a subgaussian tail estimate
$$
\Cal P \, \{ V_w - V_{\bar w} > a \} \le 
\exp \left ( -c \frac{a^2}{\tilde d^2 (w, \bar w)} \right ),
$$
where
$$
\tilde d (w, \bar w)= \left ( \sumi \Big ( w^2(i) - \bar w^2(i) \Big )^2
\right )^{1/2}.
$$
We shall estimate the metric $\tilde d$ by a quasimetric, which is simpler 
to control.
$$
\frac{1}{\sqrt{2}} \tilde d (w, \bar w) \le d (w, \bar w) = 
\left ( \sumi \Big ( w(i) - \bar w(i) \Big )^2 \cdot 
\Big ( w^2(i) + \bar w^2(i) \Big )
\right )^{1/2}.
$$
Since
$$
\align
 d (w, \bar w) =&
\left ( \sumi \frac12 \Big ( w(i) - \bar w(i) \Big )^2 \cdot 
\Big ( (w(i)+ \bar w(i))^2  + (w(i) - \bar w(i))^2 \Big )
\right )^{1/2} \le \\
& \frac{1}{\sqrt{2}} \cdot \left (\tilde d (w, \bar w) +
\nor{w- \bar w}_{\ell_4^M}^2 \right ) \le \sqrt{2} \cdot d (w, \bar w),
\endalign
$$
we have a generalized triangle inequality for $d$.
Namely for all $u,w, \bar w \in W$ 
$$
d (w, \bar w) \le 4 \cdot ( d (w, u) +  d (u, \bar w) ). \tag 3.2
$$
The balls in the quasimetric $d$ are not convex. 
However, we have the following 
\proclaim{Lemma 3}
For all $w \in W$ and $\rho >0$ 
$$
\text{\rm conv } B_{\rho}(w) \subset B_{4 \rho}(w).
$$
\endproclaim
Here we denote by $B_{\rho}(w)$ a $\rho$-ball in the quasimetric $d$.
\demo{Proof}
Note that since for all $u \in B_{\rho}(w)$
$$
\align
\left ( \sumi \Big ( u(i) - w(i) \Big )^2 w^2(i) \right )^{1/2} & \le \rho \\
\intertext{and}
\left ( \sumi \Big ( u(i) - w(i) \Big )^4  \right )^{1/4} & \le 
(\sqrt{2} \rho )^{1/2},
\endalign
$$
the same inequalities hold also for all $u \in \conv B_{\rho}(w)$.
Since for all $a,b \in \Bbb R, \ a^2+b^2 \le 4a^2 +2(a-b)^2$,
for any $u \in \conv B_{\rho}(w)$  we have
$$
\align
d(u,w) \le  & \left ( \sumi \Big ( u(i) -w(i) \Big )^2 \cdot 
\Big ( 4w^2(i)+ 2(u(i)-w(i))^2 \Big ) \right )^{1/2} \le \\
2 \cdot & \left ( \sumi \Big ( u(i)-w(i) \Big )^2 \cdot w^2(i) \right )^{1/2}
+ \sqrt{2} \cdot \left ( \sumi \Big ( u(i) -w(i) \Big )^4 \right )^{1/2}
\le 4 \rho. 
\endalign
$$
\hfil \qed

Denote 
$$
Q = \nor{P_W: \ell_1^M \to \ell_2^M}.
$$
Let now $r$ be a natural number to be chosen later. 
Let $k_0$ and $k_1$ be the largest  natural numbers so that
$$
\align
&r^{-k_0} \ge \text{diam \,} (W_1,\nor{\cdot}_{\infty}) = Q \\
&r^{-k_1} \ge \frac{Q}{\sqrt{n}}.
\endalign
$$
Then $k_1-k_0  \le (2 \log r)^{-1} \log n$.

Define functions $\varphi_k: W_1 \to \Bbb R$ by
$$
\alignat 2
\varphi_k(w)=& \min \{ \nor{u}^2 \ \Big | \ u \in \conv B_{2 r^{-k}}(w) \} +
\frac{k-k_0}{\log M}, 
&& \qquad \text{if } k=k_0 \etc k_1, \\
\varphi_{k}(w)=& 1+ \frac{1}{2 \log r} + \sum_{l=k_1}^k r^{-l} \cdot 
\frac{\sqrt{n \cdot \log ( 1+2 \sqrt{2} r^l  )}}
{Q \cdot \sqrt{\log M}},
&& \qquad \text{if } k>k_1.
\endalignat
$$
For any $w \in W_1$ the sequence $\{ \varphi_k(w) \}_{k=k_0}^{\infty} \,$
is nonnegative nondecreasing and bounded by an absolute constant 
depending only on $r$.
Indeed, if $k \le k_1$ then 
$$
\varphi_k(w) \le 1+ \frac{1}{2 \log r}\cdot \frac{\log n}{\log M}.
$$
For $k>k_1$ we have 
$$
\align
\varphi_k(w) &\le  1+ \frac{1}{2 \log r} +\sum_{l=k_1}^{\infty} r^{-l} \cdot 
\frac{\sqrt{n \cdot \log ( 1+ 2 \sqrt{2} r^l )}}{Q \cdot \sqrt{\log M}} \le \\
& 1+ \frac{1}{2 \log r} + c(r) \cdot r^{-k_1} \cdot \frac{\sqrt{n}}{Q} \cdot
\frac{\sqrt{\log ( 1+2 \sqrt{2} r^{k_1} )}} {\sqrt{\log M}} \le C(r).
\endalign
$$

To prove Lemma 1 we have to show that  condition (3.1) holds for
$\{ \varphi_k(w) \}_{k=k_0}^{\infty}$ with 
$\sigma= (c \cdot Q \cdot \sqrt{\log M})^{-1}$.
Let $x \in W_1$ and suppose that the points 
$x_1 \etc x_N \in B_{r^{-k}}(x)$ satisfy
$$
d(x_j,x_l) \ge r^{-k-1} \qquad \text{for all } j \ne l.
$$
For $k \ge k_1-1$ condition (3.1) follows from the simple volume estimate 
$$
\align
N \le & N(W_1,d,r^{-k-1}) \le 
N(W_1,\nor{\cdot}_{\infty},\frac{r^{-k-1}}{\sqrt{2}}) 
\le N(W_1, \nor{\cdot},\frac{r^{-k-1}}{\sqrt{2}} ) \le \\
& \left ( 1+ \frac{2 \sqrt{2}}{r^{-k-1}} \right )^n.
\endalign
$$

Suppose now that $k_0 \le k < k_1-1$.
For $j= 1 \etc N$ denote by $z_j$ the point of $\conv B_{2 r^{-k-2}}(x_j)$
for which the minimum of $\nor{z}$ is attained and denote by $u$ the 
similar point of $\conv B_{2 r^{-k}}(x)$.
By (3.2) and Lemma 3 we have for all $j \ne l$
$$ 
d(x_j,x_l) \le 16 \cdot \Big ( d(x_j,z_j) + d(z_j, z_l) +d(z_l,x_l) \Big )
\le 16 \cdot \Big ( 16 \cdot r^{-k-2} + d(z_j, z_l)\Big ),
$$
so, $d(z_j,z_l) \ge \frac12 r^{-k-1}$ if $r \ge 512$.  
Under the same assumption on $r$ we have 
$$
d(z_j, x) \le 4 \Big ( d(z_j, x_j) + d(x_j, x) \Big ) \le 2 r^{-k}.
$$
Denote 
$$
\th =\max_{j=1\etc N} \nor{z_j}^2 - \nor{u}^2.
$$
We have to prove that
$$
r^{-k} \cdot \left ( c \cdot Q \cdot \sqrt{\log M} \right )^{-1} \cdot
\sqrt{\log N}
\le \max_{j=1 \etc N} \varphi_{k+2}(x_j)-\varphi_k(x)= \th+ \frac{2}{\log M}.
\tag 3.3
$$

Since $\frac{z_j+u}{2} \in \conv B_{2 r^{-k}}(x)$ and $\nor{u} \le \nor{z_j}$,
we have
$$
\nor{\frac{z_j-u}{2}}^2 = \frac12 \nor{z_j}^2 + \frac12 \nor{u}^2 -
\nor{\frac{z_j+u}{2}}^2 \le \nor{z_j}^2 - \nor{\frac{z_j+u}{2}}^2 \le
\nor{z_j}^2 - \nor{u}^2,
$$
so,
$$
\nor{z_j-u} \le 2 \sqrt{\th}. \tag 3.4
$$
Thus, $N$ is bounded by the $\frac12 r^{-k-1}$-entropy of the set 
$K=u+2\sqrt{\th} \BM \cap W$ in the quasimetric $d$.
To estimate this entropy we partition the set $K$ into $S$ disjoint subsets
having diameter less than $\frac{1}{16} r^{-k-1} \th^{-1/2}$ in the 
$\ell_{\infty}$ metric.
By part (1) of Lemma 2 we may assume that
$$
\frac{1}{16} r^{-k-1}  \cdot \th^{-1/2} \sqrt{\log S} \le 
c \cdot Q \cdot \sqrt{\th} \sqrt{\log M}. \tag 3.5
$$
If $S \ge \sqrt{N}$, we are done, because in this case (3.5) implies (3.3).
Suppose that  $S \le \sqrt{N}$.
Then there exists an element of the partition containing at least $\sqrt{N}$
points $z_j$.
Let $J \subset \{ 1 \etc N \}$ be the set of the indices of these points.
We have 
$$
\nor{z_j-z_l}_{\infty} \le \frac{1}{16} r^{-k-1} \cdot \th^{-1/2}  \tag 3.6
$$
for all $j,l \in J, \ j \ne l$.
Since $d(z_j,z_l) \ge \frac12 r^{-k-1}$, we have
$$
\aligned
 \left ( \frac12 r^{-k-1} \right )^2 \le &\sumi  \Big ( z_j(i)-z_l(i) \Big )^2
\cdot \Big ( z_j^2(i)+z_l^2(i) \Big ) \le \\
& \sumi  \Big ( z_j(i)-z_l(i) \Big )^2 \cdot \\
 \Big [ 4u^2(i) +  
z_j^2(i) \cdot &\bold 1_{\{i \bigm | |z_j | \ge 2 |u(i)| \}}(i) + 
z_l^2(i) \cdot \bold 1_{\{i  \bigm | |z_l | \ge 2 |u(i)| \}}(i) \ \Big ] . 
\endaligned \tag 3.7
$$
Then (3.4) implies 
$$
\sumi z_j^2(i) \cdot \bold 1_{\{i \bigm | |z_j | \ge 2 |u(i)| \}}(i)
\le 16 \th \tag 3.8
$$
Combining (3.6) and (3.8) we get that (3.7) is bounded by
$$
2 \cdot 16 \th \cdot \left ( \frac{\th^{-1/2}}{8} r^{-k-1} \right )^2 +
4 \sumi \Big ( z_j(i)-z_l(i) \Big )^2 \cdot u^2(i).
$$
Thus, for all $j,l  \in J, \ j \ne l$ we have 
$$
\left ( \sumi \Big ( z_j(i)-z_l(i) \Big )^2 \cdot u^2(i) \right )^{1/2}
\ge \frac18 r^{-k-1}.
$$
Then part (2) of Lemma 2 implies
$$
\frac18 r^{-k-1} \sqrt{\log |J|} \le C \sqrt{\th} \cdot Q \cdot
\left ( \sumi u^2(i) \right )^{1/2} \le  C \sqrt{\th} \cdot Q.
$$
Since for all $\th>0$
$$
2 \sqrt{\th} \le \sqrt{\log M} \cdot \th +\frac{1}{\sqrt{\log M}},
$$
we get
$$
\frac{1}{16} r^{-k-1} \sqrt{\log N} \le
\frac{1}{8} r^{-k-1} \sqrt{\log |J|} \le  C  \cdot Q \cdot \sqrt{\log M}
\cdot \left ( \th +\frac{1}{\log M} \right ). \qed
$$
\enddemo
\enddemo

\subhead   Acknowledgment  \endsubhead
I would like to thank Joram Lindenstrauss  for helpful
discussions and Michel Talagrand for the permission to present his 
approach to the proof 
of Lemma~1.

\enddocument